\documentclass[a4paper, 10pt, reqno]{amsart}

\usepackage{ amssymb, amsmath, amsthm}
\usepackage{graphicx, psfrag, setspace, subfigure, gensymb}
\usepackage{amsfonts}
\usepackage{bbm}
\usepackage{fix-cm}
\usepackage{comment}
\usepackage{tikz,tikz-cd}
\usetikzlibrary{matrix,arrows,decorations.pathmorphing,decorations.pathreplacing}
\tikzset{commutative diagrams/diagrams={baseline=-2.5pt},commutative diagrams/arrow style=tikz}
\usepackage[colorlinks]{hyperref}
\usepackage{float}
\usepackage{setspace} 

\newcommand\Z{\mathbb Z}
\newcommand\C{\mathbb C}

\newcommand{\cU}{\mathcal{U}}
\newcommand{\cV}{\mathcal{V}}

\newcommand{\cX}{\mathcal{X}}

\newcommand\isoto{\stackrel{\sim}{\To}}
\newcommand\id{\mathrm 1}

\newcommand\To{\longrightarrow}

\newcommand\Hom{\operatorname{Hom}}

\renewcommand\P{\mathbb P}

\newcommand\GL{\operatorname{GL}}

\newcommand{\Sym}{\operatorname{Sym}}

\newcommand{\Wedge}{\mbox{\scalebox{1.2}{$\wedge$}}}

\newcommand{\sslash}{/\!/}

\newcommand{\al}[1]{\begin{align*}#1\end{align*}}

\newcommand{\beq}[1]{\begin{equation}\label{#1} }
\newcommand{\eeq}{\end{equation}}
\newcommand{\pgap}{\vspace{5pt}}

\theoremstyle{plain}

\newtheorem{thm}[equation]{Theorem}
\newtheorem{lem}[equation]{Lemma}

\theoremstyle{remark}
\newtheorem{rem}[equation]{Remark}

\theoremstyle{definition}

\newtheorem{eg}[equation]{Example}

\makeatletter \@addtoreset{equation}{section} \makeatother

\let\oldtocsection=\tocsection
\let\oldtocsubsection=\tocsubsection
\let\oldtocsubsubsection=\tocsubsubsection
\renewcommand{\tocsection}[3]{\hspace{0em}\oldtocsection{#1}{#2}{#3}}
\renewcommand{\tocsubsection}[3]{ \hspace{1em} \oldtocsubsection{#1}{\small{#2}}{\small{#3}} }
\renewcommand{\tocsubsubsection}[3]{\hspace{2em}\oldtocsubsubsection{#1}{\small{#2}}{\small{#3}}}

\newcommand{\marginparstretch}{0.6}
\let\oldmarginpar\marginpar
\renewcommand\marginpar[1]{\-\oldmarginpar[\framebox{\setstretch{\marginparstretch}\begin{minipage}{\marginparwidth}{\raggedleft\scriptsize #1}\end{minipage}}]{\framebox{\setstretch{\marginparstretch}\begin{minipage}{\marginparwidth}{\raggedright\scriptsize #1}\end{minipage}}}}

\newcommand{\aand}{\quad\quad\mbox{and}\quad\quad}
\newcommand\ie{\emph{i.e.}~}

\newcommand\Rep{\operatorname{Rep}}

\begin{document}

\title{The McKay correspondence in type $D_4$ via VGIT}
\author{Tarig Abdelgadir and Ed Segal}

\maketitle

\begin{abstract}We present an explicit GIT construction which produces both the minimal resolution of the type $D_4$ surface singularity, and also the orbifold resolution. Our construction is based on a Tannakian approach which is in principle applicable to arbitrary quotient singularities. 
\end{abstract}

\tableofcontents

\section{Introduction}

Let $\Gamma$ be a finite subgroup of $SL_2(\C)$ and let $X=\C^2/\Gamma$ be the corresponding Kleinian surface singularity. The McKay correspondence famously relates the representation theory of $\Gamma$ to the geometry of the minimal resolution $\widetilde{X}\to X$.  Initially a combinatorial observation \cite{McKay}, it is now understood that we should consider the orbifold 
$$ \cX = [\C^2/\Gamma] \To X $$
as another possible crepant resolution of $X$, and then other formulations of the correspondence should follow. For example there is an equivalence of derived categories $D^b(\widetilde{X}) \cong D^b(\cX)$ \cite{BKR}, which is is an instance of the general prediction that all crepant resolutions of a given singularity should be derived equivalent. There are also results relating Gromov-Witten invariants of $X$ and of $\cX$ \cite{Ruan, BA}. 
\pgap

A more elementary question we can ask is whether this birational equivalence
$$X \dashleftarrow\dashrightarrow \cX$$
can be realised via variation-of-GIT. This means looking for some larger variety $Z$,  with an action of a reductive group $G$, such that both $X$ and $\cX$ appear as possible GIT quotients $Z\sslash_\theta G$ for different choices of stability condition $\theta$. 
\pgap

For type $A$ this question is easy. All the spaces involved are toric, and standard toric geometry techniques produce the requested construction immediately, with $G$ a torus and $Z$ a vector space (see Example \ref{eg.An}). But for types D and $E$ where $\Gamma$ is non-abelian it is much less obvious how to proceed.
\pgap

In this paper we write down an explicit $Z$ and $G$ which satisfy our question for the simplest non-abelian case, $D_4$. Here $\Gamma$ is the quaternion group $Q$, \emph{aka} the binary dihedral group $BD_2$. It sits in $SL_2$ as the double cover of the dihedral group $D_2 \cong (\Z_2)^2\subset SO_3$:
$$\begin{tikzcd} Q \ar[hookrightarrow]{r} \ar{d}& SL_2 \ar{d} \\
(\Z_2)^2 \ar[hookrightarrow]{r} & SO_3
\end{tikzcd}$$
Our construction of $Z$, presented in Section \ref{sec.construction}, is a little involved. It is the vanishing locus of an explicit set of equations in a 13-dimensional vector space, carrying a natural action of the group $G=(\C^*)^3\times GL_2$. Moreover our proof that both $X$ and $\cX$ do arise as GIT quotients of $Z$ by $G$ requires quite of lot of hands-on verification.

But behind the messy details there lies a general ansatz, outlined by the first author and Chan in \cite{AC} and reviewed in our Section \ref{sec.Tannakian}. And this ansatz applies in principle to \emph{any} finite quotient singularity. So it would be a worthwhile goal to replace our hands-on proof with a more abstract and elegant one, which might then work in greater generality. 

\begin{rem}
It is perhaps surprising, given the extensive literature on the McKay correspondence, that such a construction has not appeared before. The reason is that most previous approaches have focused on constructing the minimal resolution $X$ as a moduli space of \emph{modules} (or sheaves, or quiver representations, or $G$-clusters). And since modules never have finite non-trivial stabilizer groups this approach can never produce the orbifold resolution.
\end{rem}

\begin{rem} In this paper when we discuss `GIT quotients' what we really mean is the stack-theoretic quotient of the semistable locus. Over the stable locus this just means taking the quotient as a DM stack rather than as a variety with finite quotient singularities, which is obviously essential for our purposes.

If there are strictly semistable points then the quotient stack differs more radically from the quotient variety, since it is an Artin stack, and S-equivalent orbits are not identified. But we shall show that are no strictly semistable points in our construction.
\end{rem}

\begin{rem} Ideally our construction would have brought the type $D_4$ McKay correspondence into the framework of \textit{Gauged Linear Sigma Models}, with all the accompanying string-theoretic techniques. Unfortunately, since our $Z$ is not a complete intersection (see Remarks \ref{notlci} and \ref{alleqs}) it's not clear that we have achieved this. 
\end{rem}

\subsection*{Acknowledgements}
The first author would like to thank Daniel Chan for enlightening discussions relating to this article.
This article was initiated at the ICTP, the authors would like to thank them for the hospitable enviroment they provided.

\section{The Tannakian approach}\label{sec.Tannakian}

Suppose we have a finite group $\Gamma$ acting linearly on a vector space $V$, and a corresponding orbifold $\cX=[V/\Gamma]$. How can we produce $\cX$ as a GIT quotient?

 As stated this question is trivial since $\cX$ is the GIT quotient of $V$ by $\Gamma$. But what if we ask for a construction in which $\cX$ is one of several possible quotients? Or we can simplify by ignoring $V$. Then the question becomes something like: how can we produce the stack $B\Gamma$ as a GIT quotient $Z\sslash_\theta G$, where the group $G$ is some reductive group with a character lattice of positive rank?
\pgap

Tannakian duality gives us a heuristic way to approach this problem. What we need is a presentation of the representation category $\Rep\Gamma$. By this we mean something like the following data:

\begin{enumerate}
\item A finite list of irreps $U_1,..., U_k$ which generate all $\Gamma$-representations under tensor products and Schur powers.

\item A finite list of isomorphisms
$$B_t : \cU_t \isoto \cV_t $$
where $\cU_t$ and $\cV_t$ are some given combinations of the $U_i's$, formed by sums, products and Schur powers. These $B_i$'s should generate all isomorphisms in $\Rep \Gamma$.

\item A finite list of relations that hold between the $B_t$'s, generating all such relations.

\end{enumerate}

From this data we can build a GIT problem, as follows. Treat the $U_i$'s just as a list of vector spaces, of dimensions $d_1,..., d_k$. Then the $B_t$'s collectively form an element of the vector space
$$H =  \bigoplus_{t=1}^s \Hom(\cU_t, \cV_t) $$
which carries a natural action of the group
$$G = GL_{d_1}\times ... \GL_{d_k} $$
and the relations (3) cut out a $G$-invariant subvariety $Z\subset H$. There is an open set $H^o\subset H$ where each $B_t$ is an isomorphism and a corresponding open set:
$$Z^o=Z\cap H^o$$
Tannakian duality guarantees that the stabilizer group of any point in $Z^o$ is our finite group $\Gamma$, so at the very least we have a point in the GIT quotient stack that has the correct isotropy. What we'd still have to check is that:
\begin{itemize}
\item $Z^o$ is the semistable locus for some stability condition $\theta$, and
\item $G$ acts transitively on $Z^o$. 
\end{itemize}
If these two conditions hold then the GIT quotient $Z\sslash_\theta G$ is $B\Gamma$. 
\pgap

Let's suppose that we have solved this problem, and return to our original problem of constructing $[V/\Gamma]$. The vector space $V$ is a $\Gamma$-representation, so some combination of the irreps $U_1,..., U_k$, and there is a corresponding $G$-representation which we could also denote by $V$. The obvious candidate solution to our problem is to take the variety $Z\times V$  with its action of $G$. If we can find a stability condition $\theta$ such that the $\theta$-semistable locus is $Z^o\times V$, and $G$ acts transitively on $Z^o$, then: 
$$Z\sslash_\theta G = [V/\Gamma]$$
And there are certainly other possible GIT quotients since the the character lattice of $G$ has rank $k$. We can hope to find geometric resolutions of $V/\Gamma$ amongst these other possible quotients. 

\begin{rem} If $\Gamma$ lies in $SL(V)$ then the singularity $V/\Gamma$ is Gorenstein and the orbifold $[V/\Gamma]$ is a crepant resolution, so we could hope to get \emph{crepant} geometric resolutions amongst the other resolutions. In our main example $\Gamma=Q$ (Section \ref{sec.construction}) this does indeed happen. But we have not been able to find an \emph{a priori} reason why these other quotients - even if they are resolutions - should be crepant.
  \end{rem}

\begin{eg}\label{eg.An} Let $\Gamma$ be the cyclic group $C_n$. Then $\Rep \Gamma$ is generated by a single 1-dimensional irrep $U$, subject to the single isomorphism:
$$B: U^{\otimes n} \isoto \C $$
There are no relations to impose on $B$. 

Now we can build our GIT problem. The group $G$ is $\GL(U)\cong \C^*$, acting on the vector space $H=\Hom(U^{\otimes n}, \C)$. So $H$ is one-dimensional with a $\C^*$ action of weight $-n$, and $Z=H$. For one of the two possible stability conditions it is indeed true that the semistable locus is $Z^0\cong \C^*$, and the $G$-action on this is transitive, so 
$$Z\sslash G\cong BC_n$$
(the other GIT quotient is empty). 

Now consider the type $A_{n-1}$ Kleinian singularity $X = V/C_n$, where: $$V=U\oplus U^{\otimes n-1}$$
 Our candidate for constructing resolution of this singularity is the larger GIT problem $Z\times V$, which is simply $\C^3$ carring a $\C^*$ action of weights $(-n, 1, n-1)$. We see that one of the GIT quotients is indeed the orbifold resolution $[V/C_n]$. 

The other generic GIT quotient is a non-affine orbifold given by the total space of the canonical bundle on weighted $\P^1_{1:n-1}$. This is another crepant resolution of the singularity $X$. 
\pgap

If we want the fully geometric (\ie non-orbifold) resolution we need a more redundant presentation of $\Rep C_n$. We take all the non-trivial characters $U_1,..., U_{n-1}$, and isomorphisms:
$$B_i: U_1\otimes U_i \isoto U_{i+1}\;\; \mbox{ for } i\in [1, n-2],\quad\quad B_{n-1}: U_1\otimes U_{n-1}\isoto \C$$
There are still no relations. This leads us to a GIT problem 
$$Z \times V \cong \C^{n+1} $$
with an action of the torus $\prod_i GL(U_i)\cong (\C^*)^{n-1}$. All GIT quotients will be toric surfaces, and indeed this is just the standard  construction of the surface $\widetilde{X}$ which minimally resolves the $A_n$ singularity. The orbifold $[V/C_n]$ appears as another of the possible quotients.
\end{eg}
\vspace{10pt}

In the previous example we avoided having relations on the $B_t$'s. But when we move to non-abelian cases the relations are essential, and in practice are the hardest part of implementing our ansatz.

\begin{eg} \label{eg.S3} Let $\Gamma= S_3$. We present this case both because it is the smallest non-abelian group but also because it shares some similarities with the group $Q$ which is the focus of Section \ref{sec.construction}. 

The irreps of $S_3$ can be generated from the single 2d irrep $U$, and there is an isomorphism:
$$B: \Sym^2 U \To U \oplus \C $$
The second component of $B$ gives an inner product 
$$b: U \isoto U^\vee$$
which extends to an inner product on $U\oplus \C$. But the space $\Sym^2 U$ also carries something close to an inner product, namely the canonical isomorphism:
$$J: \Sym^2 U \isoto \Sym^2 U^\vee \otimes (\det U)^2 $$
The map $B$ is, roughly speaking, an isometry with respect to these inner products. The precise statement is that the following relation holds:
\beq{S3rel}B^\vee\circ (b, 1) \circ B = (\det B) J \eeq
Note that  $\det B\in (\det U)^{-2}$. 

We believe that the data of $U, B$ and the relation \eqref{S3rel} is a presentation of the category $S_3$-rep. If we accept this claim, our ansatz leads us to consider the vector space
$$  H = \Hom(\Sym^2 U, U\oplus \C) $$
with its action of $G=GL(U)\cong GL_2$, and the $G$-invariant subvariety $Z\subset H$
cut out by \eqref{S3rel}. It's not hard to check directly that $G$ does indeed act transitively on the open set $Z^o$, with quotient:
$$[Z^o/G]=BS_3$$
But we have not attempted the stability analysis. 
\pgap

If one wanted to pursue this example further it would be interesting to set $V=U\oplus U$ since the quotient $V/S_3$ is (up to a trivial factor) the same as the 3rd symmetric product of $\C^2$. Thus one might hope to realise the birational transformation
$$[(\C^2)^{\oplus 3} / S_3 ] \dashleftarrow\dashrightarrow \mbox{Hilb}_3(\C^2) $$
via VGIT.
\end{eg}

\section{The construction}\label{sec.construction}

In this section we apply the heuristics of Section \ref{sec.Tannakian} to the quaternion group $Q$ and the corresponding $D_4$ Kleinian singularity $V/Q$. 
\pgap

First we note that $Q$ has four non-trivial irreps of dimensions 1,1,1 and 2. So we take four vector spaces $L_1, L_2, L_3$ and $V$, where $\dim L_i=1$ for each $i$ and $\dim V=2$. Then - guided by isomorphisms that hold in the representation ring - we form the vector space
$$ H = \{ (\alpha_1,\alpha_2, \alpha_3, \beta, B)\} $$
where:
\begin{align*}&\alpha_i \in \Hom(L_i^{2},\, \det V) \\
&\beta\in \Hom\big((\det V)^2,\, L_1 L_2 L_3 \big)\\
&B \in \Hom\big( \Sym^2 V, \; \bigoplus_i L_i \big) \end{align*}
Note that $\dim H = 1+1+1+1+9 = 13$ and it carries a natural action of the group:
\al{G &= \GL(L_1)\times \GL(L_2)\times \GL(L_3)\times \GL(V) \\
& \cong (\C^*)^3\times \GL_2}

Next we need to identify the relations that hold between these isomorphisms in the category $\Rep Q$, so we can specify our invariant subvariety $Z\subset H$. For this we need to introduce a little more notation.
\pgap

Let $J$ denote the canonical isomorphism
$$J: \Sym^2 V \isoto \Sym^2 V^\vee \otimes (\det V)^{\otimes 2} $$
as in Example \ref{eg.S3}. We observe that the $\alpha_i$'s give us a similar structure on the space $\bigoplus_i L_i$, since we can assemble them into a linear map:
$$A =\begin{pmatrix}
\alpha_1 & 0 & 0 \\
0 & \alpha_2 & 0 \\
0 & 0 & \alpha_3 \end{pmatrix} \!:\;  \bigoplus_i L_i \To  \bigoplus_i L_i^\vee\otimes\det V   $$
\pgap

The equations we want to write are, approximately, the statement that $B$ is an isometry with respect to these two inner products. In fact in the open set
\beq{Ho}H^o =\{\alpha_i\neq 0\,  \forall i,\; \beta\neq 0,\; \det B \neq 0\}\; \subset H\eeq
that is exactly the condition we want. But to extend over the whole of $H$ we need the following three equations, which are the key ingredient in our construction.
\pgap

\fbox{
 \addtolength{\linewidth}{-10\fboxsep}%
 \addtolength{\linewidth}{-3\fboxrule}%
 \begin{minipage}{\linewidth}
\vspace{7pt}
\beq{eq1} B^\vee A B = \alpha_1\alpha_2\alpha_3\beta^2 J\tag{E1}\eeq

\beq{eq2} BJ^{-1}B^\vee A  = \alpha_1\alpha_2\alpha_3\beta^2 I  \tag{E2}\eeq

\beq{eq3} \Wedge^2 B = \beta A B J^{-1}\tag{E3}\eeq
\vspace{1pt}
\end{minipage}
}

\vspace{10pt}

In the second equation $I$ denotes the identity map on $\oplus_i L_i$.  Implicit in third equation is the canonical isomorphism between $\Wedge^2 (\Sym^2 V)$ and $\Sym^2 V^\vee\otimes (\det V)^3$,  which means we can read $\Wedge^2 B$ as a map:
$$\Wedge^2 B: \; \Sym^2 V^\vee \To \big(L_2L_3\oplus L_1L_3\oplus L_1L_2\big)\otimes (\det V)^{-3} $$
\begin{rem}
We can write our equations more concisely by denoting
$$\omega = \alpha_1\alpha_2\alpha_3\beta^2 \aand C = J^{-1}B^\vee A$$
so \ref{eq1}, \ref{eq2}, \ref{eq3} become simply:
$$CB = \omega, \quad\quad\quad BC = \omega, \quad\quad\quad \Wedge^2 B= \beta C^\vee $$
\end{rem}
\begin{rem}\label{notlci} The subvariety $Z$ is unfortunately quite far from being a complete intersection. Indeed, in the open set $H^o$ \eqref{Ho} taking the determinant of \ref{eq3} shows that:
$$\det B = \beta^3\alpha_1\alpha_2\alpha_3 $$
But if $B$ is invertible then 
$$\Wedge^2 B = (\det B)(B^\vee)^{-1}$$
automatically, so \ref{eq1} and \ref{eq2} follow immediately from \ref{eq3}. 

However, outside $H^o$ the first two equations are independent of the third, and our construction really does require all of them. See Remark \ref{alleqs} for a full justification of this claim.
\end{rem}
\pgap

Our VGIT construction is the affine variety $Z\times V$, with the action of the group $G$. Note that a character of $G$ must be of the form $L_1^{\theta_1}L_2^{\theta_2}L_3^{\theta_3}(\det V)^{\theta_4}$ so is specified by four integers. The main result of this paper is:

\begin{thm}\label{thm.main} Let $\vartheta$ be the character $(1,1,1,1)$. 
\begin{enumerate}\setlength{\itemsep}{5pt}
\item The GIT quotient $(Z\times V) \sslash_{-\vartheta}\, G$ is the orbifold $[V/ Q]$. 
\item  The GIT quotient $(Z\times V) \sslash_\vartheta\, G$ is the minimal resolution of $V/ Q$. 
\end{enumerate}
\end{thm}

We split the proof of this theorem into the following five lemmas.
\pgap

Let $Z^o\subset Z$ be the intersection of $Z$ with $H^o$. This is the locus in $Z$ where $B$ is invertible or equivalently where $\alpha_1\alpha_2\alpha_3\beta\neq 0$. 

\begin{lem}\label{BQ} The stack $[Z^o/G]$ is equivalent to $BQ$. 
\end{lem}
Of course this is exactly what our Tannakian approach was supposed to achieve, but we need to check it since we haven't proven that our chosen data really do give a presentation of $\operatorname{Rep}(Q)$.
\begin{proof} 
The subgroup $GL_2\subset G$ acts on $Z^o$ via the usual homomorphism $GL_2 \to SO_3\rtimes \C^*$. Since this is a surjection it's clear that the action of $G$ is transitive (note that in $Z^o$ the value of $\beta$ is determined from the other variables by \ref{eq3}).

Let's examine the isotropy in the quotient group:
$$ G/\Z_2\cong (\C^*)^3\times  (SO_3\rtimes \C^*) $$
The subgroup fixing the values of $(\alpha_1, \alpha_2, \alpha_3, B)$ is the group $(\Z_2)^3$, embedded diagonally in $(\C^*)^3\times O_3$, and to also fix $\beta$ we must lie in the index-two subgroup $(\Z_2)^2\subset SO_3$. Hence the isotropy group in $G$ is the double cover of this in $GL_2$, which is $Q$.
\end{proof}

\begin{lem}\label{-thetass} The $(-\vartheta)$-stable locus in $Z\times V$ is $Z^o\times V$ and there are no strictly semistable points. 
\end{lem}
\begin{proof}
A point with $\alpha_i=0$ is destablized by the action of the 1-parameter subgroup $GL(L_i)$. Now suppose that all the $\alpha_i$'s are non-zero but $\beta=0$. Then \ref{eq2} implies that the image of $B^\vee$ is an isotropic line in $\Sym^2 V^\vee$, \ie it lies in a line spanned by some degenerate quadratic form. It follows that we can find a 1-parameter subgroup of $GL(V)$ which fixes $B$ and destablizes the point.\footnote{\emph{C.f.}~the subgroup $\mu$ in Table \ref{weights}.} So all points outside $Z^o\times V$ are unstable.

The function
$$\alpha_1^2\alpha_2^2\alpha_3^2\beta^2(\det B)$$
 is $(-\vartheta)$-semi-invariant and doesn't vanish in $Z^o\times V$, so all these points are semistable. And by Lemma \ref{BQ} all these points have finite isotropy groups so there cannot be any strictly semistable points.
\end{proof}

This completes the proof of Theorem \ref{thm.main} (1). Now we turn our attention to the opposite stability condition $\vartheta$. 
\pgap

To prove part (2) of the theorem we will utilize the standard construction of the minimal resolution as a moduli space of representations of a preprojective algebra, so we take a moment to review the relevant details of that construction (see \emph{e.g.} \cite{Ginz} for more background). 
\pgap

The preprojective algebra is the basic algebra which is Morita equivalent to the skew-group ring $\Sym^\bullet V^\vee \rtimes Q$, so its representations are equivalent to sheaves on the orbifold $[V/Q]$. We consider representations of the following form:
\beq{quiver}\begin{tikzcd}[row sep=30pt, column sep = 30pt]
	& L_1 \ar[xshift=.7ex]{d}{E_1}&        \\
\C \ar[yshift=0.7ex]{r}{E_0}  &
 V  \ar[yshift=-.7ex]{l}{D_0}  
\ar[xshift=-.7ex]{u}{D_1} 
\ar[yshift=.7ex]{r}{D_2}
\ar[xshift=.7ex]{d}{D_3}
 &      L_2 \ar[yshift=-.7ex]{l}{E_2}\\
       &  L_3\ar[xshift=-.7ex]{u}{E_3}  &
\end{tikzcd}\eeq
Here $V, L_1, L_2, L_3$ are the same vector spaces of dimensions 2,1,1,1 which we've used above, and each $D_i, E_i$ is a linear map. The relations in the algebra are that $D_iE_i=0$ for each $i$, and
$$\sum_{i=1}^3 E_iD_i - E_0D_0 = 0$$
(the minus sign is an optional convention that is convenient for us). The group acting is our same group $G$, so $\vartheta$ gives a stablity condition on this moduli stack of representations. An elegant observation of King \cite{King} is that a representation is stable for this character iff it is generated from vertex 0 (the vertex decorated with $\C$), which means the following two conditions hold:
\vspace{1pt}
\begin{itemize}\setlength{\itemsep}{5pt} \item $D_i\neq 0$ for each $i=1,2,3$.
\item For some $i\in\{1,2,3\}$, $V$ is spanned by $E_0$ and the image of $E_i$. 
\end{itemize}
There are no strictly semistable representations. 
\pgap

The moduli space of stable representations turns out to be the minimal resolution $X$ of $V/Q$. Our claim is that the variety $Z$ contains enough information about the representation theory of $Q$ that we can replicate this construction. 
\pgap

First we need a little more notation. For each $i$ we have a component 
$$B_i: \Sym^2 V \to L_i$$
of $B$. We denote
$$ B_i^x = B_i(x, -) : V \to L_i $$
and recall that
 $\omega$ denotes the two-form:
$$\omega = \alpha_1\alpha_2\alpha_3\beta^2 \; \in (\det V)^{-1}$$
Then for each point $(\alpha_1,\alpha_2,\alpha_3,\beta, B,x)\in H\times V$ we associate a representation of the quiver \eqref{quiver} by setting:
\addtolength{\jot}{5pt}
\al{ E_0 =x :& \;\;\C \to V\\
D_0 =  \omega(x, -) :&\;\; V \to \C\\
D_i = B_i^x :& \; \;V \to L_i \\
E_i = \alpha_i (B_i^x)^\vee:& \;\; L_i \to V }
\pgap

Thus we have a map from $[(H\times V )/ G]$ to the stack of representations of the free quiver.

\begin{lem}
If $(\alpha_1,\alpha_2,\alpha_3,\beta, B,x)$ lies in $Z\times V$ then the associated representation satifies the preprojective algebra relations.
\end{lem}
\begin{proof} At the four external vertices the relation is essentially automatic; if $L$ is a 1-dimensional vector space and $E: L \to V$ is any linear map then the composition 
$$ E^\vee\circ E \!:\; L \To L^{-1}\det V$$
is zero. At the central vertex the required relation is:
$$\sum_i \alpha_i (B_i^x)^\vee \circ B_i^x\; -\; x\circ \omega(x,-) \;=0 $$
By construction this is a trace-free endomorphism of $V$, and the space of such endomorphisms can be canonically identified with $\Sym^2 V\otimes (\det V)^{-1}$. The relation is \ref{eq1} evaluated at the point $x^2\in \Sym^2 V$. 
\end{proof}

\begin{lem}\label{thetastability} Given a point 
$$(z,x) = (\alpha_1,\alpha_2,\alpha_3,\beta, B,x)\;\in Z\times V$$
 the following are equivalent:
\begin{enumerate}
\item[(i)] The point $(z,x)$ is $\vartheta$-semistable. 
\item[(ii)] The associated representation is $\vartheta$-stable. 
\item[(iii)] We have $B_i^x\neq 0$ for each $i$, and for some $i$ the space $V$ is spanned by $x$ and the image of $\alpha_i (B_i^x)^\vee$.
\end{enumerate}
\end{lem}
\begin{proof}
That (iii)$\Rightarrow$(ii) is immediate from our discussion of King stability above. To see that (ii)$\Rightarrow$(i) we just note that if a representation is stable then there is a $\vartheta$-semi-invariant function on the moduli stack which doesn't vanish at this point. Pulling this function back shows that $(z,x)$ is $\vartheta$-semistable. 

It remains to prove that if $(z,x)$ does not satisfy the condition (iii) then it is unstable. We prove this by some slightly tedious analysis of the action of 1-parameter subgroups in $G$. 

First note that the 1-parameter subgroup 
$$t \mapsto (t,t,t,t\id_V) \in G $$
 destabilizes points with $x=0$, so we can assume $x\neq 0$, and then extend it arbitrarily to a basis $x,y\in V$. This basis induces co-ordinates on each space $\Sym^2 V^\vee \otimes L_i$, \ie it splits each quadratic form $B_i$ into:
$$(p_i, q_i, r_i) = \big(B_i(x^2), B_i(xy), B_i(y^2)\big) $$

For brevity we write $\lambda_i = GL(L_i)$ and we let $\mu\subset G$ be the 1-parameter subgroup defined by:
$$\mu(t): x \mapsto x, \quad \mu(t): y\mapsto ty $$
These four 1-parameter subgroups all pair positively with $\vartheta$, and generate a torus which acts on $H\times V$ but fixes $x\in V$. To show that $(z,x)$ is $\vartheta$-unstable it is sufficient to show that $z$ is destabilized by some 1-parameter subgroup in this torus. For convenience we present the weights of this torus action in Table \ref{weights} together with the weights of some particular 1-parameter subgroups. Note that the three columns of entries under each $B_i$ are the weights of $(p_i, q_i, r_i)$ respectively. 

\begin{table}\renewcommand{\arraystretch}{1.2}
\begin{tabular}{c|cccccccc}
 							& $\alpha_1$ & $\alpha_2$ & $\alpha_3$ & $\beta$ & $B_1$ & $B_2$  & $B_3$  \\ \hline
 $\lambda_1$				& -2 & 0 & 0 & 1 & (1, 1, 1) & (0,0,0) & (0,0,0) \\
$\lambda_2$ 				 & 0 & -2 & 0 & 1 & (0,0,0) & (1,1,1) & (0,0,0) \\
$\lambda_3$					 & 0 & 0 & -2 & 1 & (0,0,0) & (0,0,0) & (1,1,1) \\
$\mu$   						 & 1 & 1 & 1 & -2 & (0,-1,-2) & (0,-1,-2) & (0,-1,-2)  \\
\hline
$\mu+\lambda_1$			& -1 & 1 & 1 & 0 & (1,0,-1) & (0,-1,-2) & (0,-1,-2)  \\
$\mu+\lambda_1+\lambda_2$ & -1 & -1 & 1 & 0 & (1,0,-1) & (1,0,-1) & (0,-1,-2)  \\
$2\mu+\sum \lambda_i$ 		& 0 & 0 & 0 & -1 & (1,-1,-3) & (1,-1,-3)  & (1,-1,-3)   \\
\end{tabular}
\vspace{7pt} 
\caption{}\label{weights}
\end{table}

A point is destabilized by one of these 1-parameter subgroups if all co-ordinates carrying positive weights are zero. By inspection the following subsets are unstable:
\al{ &\{\beta=0, \; B_i=0\} \quad \mbox{for any }i \\
&\{\alpha_1=0, \; \alpha_2=0, \; \alpha_3=0\}\\
&\{\alpha_2=0,\;\alpha_3=0,\; p_1=0\} \quad \mbox{\emph{(and permutations thereof)}} \\
 &\{\alpha_3=0,\;p_1=0,\; p_2=0\} \quad \mbox{\emph{(and permutations thereof)}} \\
 &\{p_1=0,\; p_2=0,\; p_3=0\} }

Now suppose that for each $i$, the image of $\alpha_i(B_i^x)^\vee$ lies in the span of $x$. This is equivalent to $\alpha_ip_i=0$ for each $i$. It follows that $(z,x)$ lies in one of the unstable subsets in the list above.

It remains to show that if $B_i^x=0$ for any $i$ then $(z,x)$ is unstable, and for concreteness let us suppose $B_1^x=0$, \ie that $p_1=q_1=0$.  We consider two cases.

 Firstly suppose that $B_1=0$, \ie that $r_1=0$ too. Then \ref{eq3} implies that
$$0 = B_1\wedge B_2 = \beta\alpha_3 B_3 J^{-1} \aand 0 = B_1\wedge B_3 = \beta\alpha_2 B_2 J^{-1} $$
so in particular
$$ \beta\alpha_2p_2 = \beta\alpha_3p_3=0$$
from which it follows again that $(z,x)$ lies in one of the unstable subsets listed above.

 Secondly, suppose that $B_1\neq 0$, \ie that $r_1\neq 0$. Then since \ref{eq2} implies that 
$$ \alpha_1B_1J^{-1}B_j = 0 \aand \alpha_jB_1J^{-1}B_j =0 \quad \quad\mbox{for }j=2, 3 $$
we must have:
$$ \alpha_1 p_2=  \alpha_1 p_3= \alpha_2p_2=\alpha_3p_3 = 0 $$
It follows again that $(z,x)$ must lie in one of the unstable subsets listed above.
\end{proof}

\begin{lem} There are no strictly $\vartheta$-semistable points in $Z\times V$, and the GIT quotient $(Z\times V)\sslash_\vartheta G$ is isomorphic to the space of $\vartheta$-stable representations of the preprojective algebra. 
\end{lem}
\begin{proof} From our previous lemmas we know we have a $G$-equivariant morphism from the semistable locus in $Z\times V$ to the space of stable representations of the algebra. We also know (Lemma \ref{thetastability}(iii)) that the semistable locus has a Zariski-open cover $\cU_1\cup\, \cU_2\cup \cU_3$, mapping to a corresponding open cover $\cV_1\cup\cV_2\cup \cV_3$ of the target. We will show by direct inspection that each map $\cU_i\to \cV_i$ is an isomorphism. This proves the lemma, in particular there cannot be any strictly $\vartheta$-semistable points since there are no strictly semistable representations and the morphism is an injection.  

 For concreteness we show that $\cU_1\to \cV_1$ is an isomorphism. In the locus $\cU_1$ we have $\alpha_1\neq 0$ and we know that $x$ and $(\alpha_1(B_1^x)^\vee\circ B_1)(x^2)$ from a basis for $V$ at all points. So if we identify $V$ with $\C^2$, and $L_1$ with $\C$, then up to a unique element of $GL(L_1)\times GL(V)$ we can set 
$$x =(1,0), \quad \alpha_1=1,\quad B_1 = (1,0,r_1)$$
for some $r_1\in \C$. So $B_1(x^2)=1$ and $\alpha_1(B_1^x)^\vee(1)=(0,1)$. If we also pick basis vectors for $L_2$ and $L_3$ then our remaining variables are
$$\alpha_2, \;\;\alpha_3, \;\;\beta,\;\; B_2 = (p_2, q_2, r_2), \;\; \mbox{and}\;\; B_3 = (p_3, q_3, r_3)$$
and the only group acting is the torus $GL(L_2)\times GL(L_3) \cong(\C^*)^2$. Semistability is the condition that $(p_2, q_2)\neq (0,0)$ and $(p_3, q_3)\neq (0,0)$. 

By \ref{eq2} we have:
\al{B_1J^{-1}B_1 &=2r_1 = \beta^2\alpha_2\alpha_3 = \omega\\
B_1J^{-1}B_2 &=r_2+r_1p_2 = 0\\
B_1J^{-1}B_2 &=r_3+r_1p_3 = 0}
so all the $r_i$ variables are redundant. And by \ref{eq3} we have
$$B_1\wedge B_2 = (-r_1q_2,\, r_1p_2 - r_2, \,q_2) = \beta\alpha_3(r_3,\, -\frac{q_3}{2},\, p_3)$$
$$B_3\wedge B_1 = (r_1q_3,\, r_3-r_1p_3, \,-q_3) = \beta\alpha_2(r_2, \,-\frac{q_2}{2},\, p_2)$$
hence:
\beq{q=p} q_2 = \beta\alpha_3p_3, \quad\quad q_3 = -\beta\alpha_2 p_2 \eeq
Finally by equation \ref{eq1} we have:
\beq{normp} 1 + \alpha_2p_2^2 + \alpha_3p_3^2 = 0 \eeq
One can verify that these equations \eqref{q=p} and \eqref{normp}, together with the equations for the $r_i$'s above, imply all the remaining components of \ref{eq1}, \ref{eq2} and \ref{eq3}. 

Now consider our open set $\cU_1$ in the space of stable representations. We can similarly fix
$$E_0 = (1,0), \quad D_1= (1,0)^\top, \quad E_1 = (0,1) $$
leaving only $(\C^*)^2$ acting. We have $D_i\neq 0$ for $i=2,3$, and we know $D_i\circ E_i = 0$. So if we write $D_i = (\hat{p}_i, \hat{q}_i)^\top$ for $i=2,3$ then we must have
$$D_0=(-\hat{\omega}, 0)^\top \aand E_i = \hat{\alpha}_i ( -\hat{q}_i, \hat{p_i}), \quad i=2,3 $$
for some $\hat{\omega}, \hat{\alpha}_2, \hat{\alpha}_3$. Moreover the relation at the central vertex implies three equations:
\al{ \hat{\alpha}_2\hat{p}_2^2 + \hat{\alpha}_3\hat{p}_3^2 + 1 & = 0 \\
\hat{\alpha}_2\hat{p}_2\hat{q}_2 + \hat{\alpha}_3\hat{p}_3\hat{q}_3 & = 0\\
\hat{\alpha}_2\hat{q}_2^2 + \hat{\alpha}_3\hat{q}_3^2 - \hat{\omega} &=0 }
The first implies that $(\hat{\alpha}_2\hat{p}_2, \hat{\alpha}_3\hat{p}_3)\neq (0,0)$, then the second implies that 
$$(\hat{q}_2, \hat{q}_3) = \hat{\beta}(\hat{\alpha}_3\hat{p}_3, -\hat{\alpha}_2\hat{p}_2) $$
for some $\hat{\beta}$, and then the third becomes equivalent to $\hat{\omega} = \hat{\beta}^2\hat{\alpha}_2\hat{\alpha}_3$. Comparing with \eqref{q=p}, \eqref{normp} above it is evident that $\cV_1$ is isomorphic to $\cU_1$. 
\end{proof}

\begin{rem}\label{alleqs} Let us conclude by justifying the necessity of all three equations \ref{eq1}, \ref{eq2} and \ref{eq3}, as promised in Remark \ref{notlci}. 
 
Firstly, equation \ref{eq3} is clearly necessary since without it the variety $Z$ would have (at least) two irreducible components related by $\beta\mapsto -\beta$. 

Now consider the subset 
$$H' = \big\{\beta=0, \alpha_i\neq 0\; \forall i\big\} \subset H$$ 
In this subset:
\begin{itemize}\item[] \ref{eq1} says that the image of $B$ is isotropic in $\bigoplus L_i$.
\item[] \ref{eq2} says that the image of $B^\vee$ is isotropic in $\Sym^2 V^\vee$.
\item[] \ref{eq3} says that $B$ has rank 1.
\end{itemize}
The first two are independent, but either of them imply the third. 

Take a generic point $(u, v)\in H'\times V$ where \ref{eq2} holds but \ref{eq1} fails. Then the associated quiver representation will not satisfy the preprojective algebra relations, but will still be generated from vertex 0 and hence King stable. So $(u,v)$ is $\vartheta$-semistable. This shows that equation \ref{eq1} is necessary for our construction. 

Now instead take a point $(u, v)$ where \ref{eq1} holds but \ref{eq2} fails. Using the functions
$$\big(B_i J^{-1} B^\vee\big)_{ij} \in L_iL_j(\det V)^{-2}$$
together with $\alpha_1,\alpha_2,\alpha_3$ it is easy to find a non-vanishing function which is semi-invariant for some power of $-\vartheta$. So $(u,v)$ is $(-\vartheta)$-semistable (strictly, since the stabilizer is infinite). Thus Lemma \ref{-thetass} would fail without equation \ref{eq2}. 
\end{rem}

\bibliographystyle{halphanum}

\end{document}